\documentclass[11pt, reqno]{amsart}
\usepackage{amssymb,amsmath,amsthm,amscd}
\input{amssym.def}

\numberwithin{equation}{section}

\newcommand{\bt}{\begin{Theorem}}
\newcommand{\et}{\end{Theorem}}
\newcommand{\bi}{\begin{itemize}}
\newcommand{\ei}{\end{itemize}}
\newcommand{\bea}{\begin{eqnarray}}
\newcommand{\eea}{\end{eqnarray}}

\theoremstyle{plain}

\newtheorem{Theorem}{\sc Theorem}[section]
\newtheorem{Lemma}{\sc Lemma}[section]
\newtheorem{Proposition}{\sc Proposition}[section]

\newtheorem{Corollary}{\sc Corollary}[section]

\theoremstyle{definition}

\theoremstyle{remark}
\newtheorem{Remark}{\sc Remark}[section]

\newcommand{\be}{\begin{equation}}
\newcommand{\ee}{\end{equation}}

\def\rep{{representation~}}%

\newcommand{\bb}[1]{\mathbb #1}

\def\C{\mbox{${\mathbb C}$}}
\def\D{\mbox{${\mathbb D}$}}

\def\R{\mbox{${\mathbb R}$}}

\newcommand{\g}{{\mathfrak g}}%

\newcommand{\inner}[2]{\langle #1,#2 \rangle }%
\newcommand{\mat}[4]{\Big ( \begin{matrix} #1 & #2\\ #3 & #4 
\end{matrix} \Big )}%

\title[classification of homogeneous operators]{A  classification of homogeneous operators\\ in the Cowen-Douglas class}

\author{Adam Kor\'{a}nyi}
\address{Lehman College\\
The City University of New York\\
Bronx, NY 10468
}
\email{adam.koranyi@lehman.cuny.edu}
\author{Gadadhar Misra}
\address{Department of Mathematics \\
Indian Institute of Science\\
Bangalore 560 012
 }
\email{gm@math.iisc.ernet.in}
\thanks{The research of the first author was supported, in part, by a DST -
NSF S\&T Cooperation Program and a PSC-CUNY grant. The research of the
second author was supported, in part, by DST and UGC through DSA - SAP - Phase IV}



\begin{document}

\begin{abstract}
An explicit construction of all the homogeneous holomorphic
Hermitian vector bundles over the unit disc $\mathbb D$ is given.
It is shown that every such 
vector bundle is a direct sum of irreducible ones.
Among these irreducible homogeneous holomorphic Hermitian vector bundles
over $\mathbb D$, the ones corresponding to operators in the Cowen-Douglas
class ${\mathrm B}_n(\mathbb D)$ are identified. The  
classification of homogeneous operators in ${\mathrm B}_n(\mathbb D)$ is 
completed using an explicit realization of these operators.
We also show how the homogeneous operators in ${\mathrm B}_n(\mathbb D)$
split into similarity classes.
\end{abstract}

\maketitle
\section{Introduction}
\noindent
An operator $T$ is said to be {\em homogeneous} if its spectrum is contained in
the closed unit disc and for every M\"{o}bius transformation $g$ of the unit
disc $\mathbb D$, the operator $g(T)$ defined via the usual holomorphic
functional calculus, is unitarily equivalent to $T$.  To every homogeneous
irreducible operator $T$ there corresponds (cf. \cite[Theorem 2.2]{shift}) an
{\em associated projective unitary representation}  $U$ of the  M\"{o}bius group
$G_0$:
$$U_g^*\, T \, U_g =  g(T),\: g \in G_0.$$
The projective unitary representations of $G_0$ lift to unitary representations of
the universal cover $\widetilde{G_0}$ which are quite well-known.  We can choose (cf.
\cite[Lemma 3.1]{shift}) $U_g$ such that $ k\mapsto U_k$ is a representation of
the rotation group. 
If
$$
\mathcal H(n) = \{ x  \in \mathcal H: U_{k_\theta} x = e^{i\,n \theta} x\},
$$
where $k_\theta(z) = e^{i\theta} z$, then $T:\mathcal H(n) \to \mathcal H(n+1)$ 
is a block shift.  A complete classification of these for $\dim \mathcal H(n) \leq 1$ 
was obtained in \cite{shift} using the representation theory of $\tilde{G}_0$.  First 
examples for $\dim \mathcal H(n) = 2$ appeared in \cite{W}.  Recently \cite{KM, KMISI}, 
an $m$ - parameter family of examples with $\dim \mathcal H(n) = m$ was constructed.  
We will use the ideas of \cite{KM, KMISI} to obtain a complete classification of 
the homogeneous operators in the Cowen-Douglas class.  
Finally, we describe the similarity classes within the homogeneous Cowen-Douglas 
operators. As a consequence, we obtain an affirmative answer to the Halmos question 
(cf. \cite{GP}) for this class of operators.  We also include a somewhat 
new conceptual presentation of the Cowen-Douglas theory and a brief description of 
the method of holomorphic induction, which will be our main tool.  Our paper is 
essentially self contained and can be read without the knowledge of \cite{KM} and 
\cite{KMISI}.  The results of this paper were announced in \cite{KMIEOT} except for
Theorem \ref{Halmos}.


\subsection{\sf Vector bundles} Let $M$ be a complex manifold and suppose $\pi:
E \to M$ is a complex vector bundle.  We write, as usual, $E_z = \pi^{-1}(z)$.
For a trivialization, $\varphi: E\to M\times \mathbb C^n$, we write $\varphi(v)
=(z, \varphi_z(v))$
for $v\in E_z$ with $\varphi_z: E_z \to \mathbb C^n$ linear.  (All we are going 
to say here would be valid using
local trivializations, but in this article we will always work with global
trivializations.)

We write $E_z^*$ for the complex anti-linear dual of $E_z$, $z\in M$,
and we write $[u,v]$ for $u(v)$, $u\in E_z^*,\, v\in E_z$.
We consider $\mathbb C^n$ to be equipped with its natural inner product and
identify it with its own anti-linear dual (so $\xi \in \mathbb C^n$ is
identified with the anti-linear map $\eta\mapsto \langle \xi , \eta
\rangle_{\mathbb C^n}$).  Then $\varphi_z^*: \mathbb C^n \to E_z^*$ is
well-defined.  We set $\psi_z = {\varphi_z^*}^{-1}$ and $\psi(u) =
(z, \psi_z(u))$ 
for $u\in E_z^*$. This makes $E^*$ a
complex vector bundle with trivialization $\psi$.  We call $\varphi$ and $\psi$,
the {\em associated trivializations} of  $E$ and $E^*$.  If $E$ is a holomorphic
vector bundle then $E^*$ is an anti-holomorphic vector bundle (meaning that for
any two trivializations, $\psi_\alpha$ and $\psi_\beta$, the transition
functions $z\mapsto (\psi_\alpha)_z \circ {(\psi_\beta)_z}^{-1}$ are
anti-holomorphic) and vice-versa.

If $E$ has a Hermitian structure, we automatically equip $E^*$ with the dual
structure (giving the dual norm of $E_z$ to $E_z^*$ for all $z \in M$).

By an automorphism of $\pi: E \to M$, we mean a  diffeomorphism
$\hat{g}: E\to E$ such that $\pi\circ \hat{g} = g\circ \pi$ for some automorphism
$g$ of $M$.  We write $g_z$
for the restriction of $\hat{g}$ to $E_z$.  The automorphism $\hat{g}$ also acts on the sections 
$f$ of $E$, by $(\hat{g}^*f)(z) = g_z^{-1}f(g z)$.  When $G$ is the group of automorphisms of 
$E$ (acting on the left, as usual) we have a representation $U$ of $G$ on the sections given by 
$U_{\hat{g}}f=(\hat{g}^{-1})^*f$, that is,
$$
(U_{\hat{g}}f)(z) =  g_z\,f (g^{-1}z).
$$
Given an automorphism $g$  of $E$, there is a corresponding automorphism
of $E^*$, where the place of $g_z$ is taken by ${g_z^*}^{-1}$.  This also
remains true in the category of Hermitian bundles.  It follows that a group $G$
of automorphisms of $E$ also acts as a group of automorphisms of $E^*$.  If $E$
is \emph{homogeneous}, that is, the action of $G$ is transitive on $M$, then so is
$E^*$, and vice-versa.

\subsection{\sf Reproducing kernel} We describe,  essentially   following
\cite{BerHil},  how the usual formalism of reproducing kernels can be adapted to
vector bundles.  Suppose $\mathcal H$ is a Hilbert space whose elements are
sections of a vector bundle $E\to M$ and suppose the maps ${\rm ev}_z\,:
\mathcal H \to E_z$ are continuous for all $z\in M$. Then setting $K_z= {\rm
ev}_z^*$, we have
\begin{equation} \label{repr}
[u, f(z)] = [u , {\rm ev}_z(f)] = \inner{K_z u}{f}_\mathcal H,\,\:\: u\in E_z,f\in
\mathcal H.
\end{equation}
For all $w\in M$ then $K_wu$ is in $\mathcal H$   and is linear in $u$.  So, we
can write $K_w(z) u = {\rm ev}_z\big (K_w u \big ) =  {\rm ev}_z {\rm
ev}_w^*(u)$. We also write $K(z,w) = K_w(z) = {\rm ev}_z{\rm ev}_w^*$ which is a
linear map $E_w^* \to E_z$, and is called the reproducing kernel of $\mathcal
H$, \eqref{repr} is the reproducing property.

Clearly, $K(w,z) = K(z,w)^*$.  We have the positivity $\sum_{j,k} [{u_k},
{K(z_k,z_j)u_j}] \geq 0$ for any $z_1, \ldots ,z_p$ in $M$ and $u_1, \ldots
,u_p \in E_z^*$ which is nothing but the inequality
$$\sum_{j,k}\inner{({\rm ev}_{z_k})^* u_k}{({\rm ev}_{z_j})^*u_j}_\mathcal H \geq 0.$$
Conversely, a $K$ with these properties is
always the reproducing kernel of a Hilbert space of sections of $E$ (cf.
\cite{BerHil}).

Suppose we have a vector bundle $E$ and a Hilbert space $\mathcal H$ of sections
of $E$ with reproducing kernel $K$; suppose $\hat{g}$ is an automorphism of $E$. Then $\hat{g}$
acts on the sections  of $E$ by $(g^*f)(z)=g_z^{-1}f(gz)$.  By the density of
linear combinations of the sections of the form $K_w u$, the condition for this action to
preserve $\mathcal H$ and act on it isometrically is
$$
\inner{g^*(K_wu)}{K_z u^\prime}_\mathcal H= \inner{K_w u} {(g^{-1})^*(K_z
u^\prime)}_\mathcal H
$$
for all $z,w;\, u,u^\prime$. Evaluating both sides using \eqref{repr}, this
amounts to
$$
K(gz, gw) = g_z\,K(z,w)\,g_w^*,\,\, \mbox{for all}\: z,\,w \in M.
$$

The following remarks will be important for us.  Suppose each ${\rm ev}_z$ is
non-singular, that is, its range is the whole of $E_z$. (This is so in the important
case where $\mathcal H$ is dense in the space of sections of $E$ in the
topology of uniform convergence on compact sets.) Then $K_z = {\rm ev}_z^*$ is
an imbedding of $E_z^*$ into $\mathcal H$. Postulating that this imbedding is an
 isometry we obtain a canonical Hermitian structure on $E^*$.  Using
\eqref{repr} we can write for the norm on $E^*$
$$
\|u\|^2_z = \|K_z u\|_\mathcal H^2 = [u,K(z,z)u],\,\:\: u \in E_z^*.
$$
The vector bundle $E$ has the dual Hermitian structure, for $v \in E_z$ we have
$$
\|v\|_z^2 = [K(z,z)^{-1}v, v].
$$
In fact this statement amounts to
$$
|[u,v]|^2 \leq [K(z,z)^{-1}v, v][u,K(z,z)u]
$$
for all $u,v$ with equality reached for some $u,v$.  Since $K(z,w)$ is bijective
by hypothesis, any $v\in E_z$  can be written as $v=K(z,z) u^\prime$ with
$u^\prime \in E_z^*$  and the inequality to be proved is equivalent to
$|[u,K(z,z)u^\prime]|^2 \leq [u^\prime, K(z,z) u^\prime][u,K(z,z)u]$. 
But this is just the Cauchy-Schwarz inequality.

When $E$ is a holomorphic vector bundle, $K(z,w)$ depends in $z$ holomorphically and on $w$ 
anti-holomorphically. Hence $K(z,w)$ is completely determined by $K(z,z)$. It follows that 
$K(z,w)$ is completely determined by the canonical Hermitian structure of $E$ (or $E^*$).

In the last paragraphs, we had a Hilbert space $\mathcal H$ of
sections of $E$ and (under the assumption that each ${\rm ev}_z$ is surjective)
we associated to it a family of imbeddings of $E_z^*$, the fibres of $E^*$, into
$\mathcal H$.  This procedure can be reversed which is of importance for what
follows. Suppose now that $E$ is a vector bundle and the fibres $E_z^*$ of $E^*$
form a  smooth family of subspaces of some Hilbert space $H$ which together span
$H$, that is, $E^*$ is a anti-holomorphic sub-bundle of the trivial bundle $M \times H$). We
write $\iota_z:E_z^* \to H$  for the (identity) imbeddings.  We define,
$\tilde{f}(z) = \iota_z^* f$ for $f \in H,\,z\in M$.  Then $\tilde{f}$ is a section of $E$
and ${\rm ev}_z(\tilde{f}) = \iota_z^*f$.  If we denote by $\mathcal H$ the Hilbert
space of all $\tilde{f}$, $f \in H$, with norm $\|\tilde{f}\| = \|f\|$, each ${\rm ev}_z$ is
continuous, so we have a reproducing kernel Hilbert space. The reproducing
kernel is $K_z u = \widetilde{\iota_z u}$.

\subsection{\sf Operators in the Cowen-Douglas class}
We modify the definition of the class of operators introduced in \cite{C-D} in
an inessential way.  A conceptual presentation in which the role of the dual of
the bundle constructed in \cite{C-D} is apparent follows.
Given a domain $\Omega \subseteq \mathbb C$, we say the
bounded operator $T$ on the Hilbert space $H$ is in ${\mathrm B}_n(\Omega)$ if
$\bar{z}$ is an eigenvalue of $T$ , the range of the operator $T-\bar{z}$ is
closed, and the corresponding eigenspaces $F_z$ are of constant dimension $n$
for $z \in \Omega$.  It is proved in \cite{C-D} that the spaces $F_z$ span an
anti-holomorphic Hermitian vector bundle $F\subseteq \Omega \times H$. (In
\cite{C-D}  the eigenvalues are $z\in \Omega$ and so $F$ is a holomorphic vector
bundle; it is more convenient for us to change this.)  We write, for $z\in
\Omega$, $\iota_z:F_z \to  H$ for the identity imbedding.  Now, $E=F^*$ is a
holomorphic vector bundle, this will be the primary object for us.  The bundle
$F$ is identified with $E^*$, in what follows we refer to it as $E^*$. We are
now in the situation discussed above in Section 1.2.

To the elements $f$ of $H$ there correspond the sections $\tilde{f}$ of $E$
(defined by $\tilde{f}(z) = \iota^*_z f$) and form a Hilbert space $\mathcal H$
isomorphic with $H$ and having a reproducing kernel
$K_z u = \widetilde{\iota_z u}$.

Under this isomorphism, the operator on $\mathcal H$ corresponding to $T$ is
$M^*$, where $M$ is the multiplication operator $(M\tilde{f})(z)  =
z\tilde{f}(z)$. In fact (cf. \cite{C-D}) for any $u \in E_z^*$, 
\begin{eqnarray*}
[u, \widetilde{T^*f}(z)] = \inner{\iota_z u}{T^*f}_H &=& \inner{T\iota_zu }{f}_H =
\bar{z}\inner{\iota_z u }{f}_H\\
&=& [u,z\tilde{f}(z)]=[u,M\tilde{f}(z)] \\
\end{eqnarray*}
\subsection{\sf Trivialization} Finally, we describe how the preceding material
appears when the vector bundle is trivialized.  We always use associated
trivializations $\varphi, \psi$ of $E$ and $E^*$.  As explained in the
beginning, this means that $\psi_z = {\varphi_z^*}^{-1}$, that is,
$[u,v]=\inner{\psi_zu}{\varphi_z v}_{\mathbb C^n}$ for $u\in E_z^*$ and $v \in
E_z$.  We will consider here only the case where $E$ is a holomorphic vector
bundle.  When $g$ is an automorphism of $E$, in the trivialization $g_z:E_z \to
E_{g z}$ becomes $\varphi_{gz} \circ g_z \circ\varphi_z^{-1}$, which we write as the
matrix $J_g(z)^{-1}$.  When $g$ is followed by another automorphism $h$, the
relation $(hg)_z = h_{gz}\circ g_z$ becomes the multiplier identity
\begin{equation}\label{multidnew}
J_{h g}(z)=J_g(z)J_h(gz).
\end{equation}
For the induced automorphism of $E^*$, the place of
$J_g(z)$ is taken by ${J_g(z)^*}^{-1}$.

The sections of $E$ (resp  $E^*$) in the trivialization become the holomorphic
(resp  anti-holomorphic)  functions $\hat{f}(z) = \varphi_z(f(z))$ (resp
$\psi_z(f(z))$).  The action $g^* f$ of an automorphism $g$ on a section becomes
$\big (g^*\hat{f}\big )(z)=J_g(z)\hat{f}(gz)$.  If $G$ is a group of automorphisms 
of $E$, the representation $U$ of $G$ described in Section 1.1 becomes the ``multiplier representation''
\begin{equation} \label{multrepnew}
\big (\hat{U}\hat{f} \big )(z) = J_{g^{-1}}(z)\hat{f}(g^{-1}z).
\end{equation}

A Hermitian structure on $E$ becomes a family of inner products on $\mathbb
C^n$, parametrized by $z\in M$.  One can always write
$$
\|\xi\|_{E_z}^2 = \inner{H(z)\xi}{\xi}_{\mathbb  C^n}
$$
with a positive definite matrix $H(z)$, $z\in M$. The structure is invariant
under a bundle automorphism $\hat{g}$ if and only if
$\inner{H(gz)J_g(z)^{-1}\xi}{J_g(z)^{-1}\xi}_{\mathbb C^n} = \inner{H(z)\xi}{\xi}_{\mathbb C^n}$, that is,
$$
H(gz) = J_g(z)^*H(z)J_g(z).
$$
The dual Hermitian structure of $E^*$ is given by $\|\xi\|_{E^*_z}^2 =
\inner{H(z)^{-1}\xi}{\xi}_{\mathbb  C^n}$.

A Hilbert space $\mathcal H$ of sections of $E$ becomes a space $\hat{\mathcal
H}$ of holomorphic functions from $M$ to $\mathbb C^n$.  The reproducing kernel
becomes $\hat{K}(z,w) = \varphi_z\circ K(z,w) \circ \psi_w^{-1}$, a matrix
valued function, holomorphic in $z$ and anti-holomorphic in $w$.  The
reproducing property appears as
$$
\inner{\hat{f}(z)}{\xi}_{\mathbb C^n}= \inner{\hat{f}}{\hat{K}_z
\xi}_{\hat{\mathcal H}},
$$
the positivity as
$$
\sum_{j,k}\inner{\hat{K}(z_j,z_k)\xi_k}{\xi_j}_{\mathbb C^n} \geq 0,
$$
and the isometry of the $G$ - action as
$$
J_g(z)\hat{K}(gz,gw)J_g(w)^* = \hat{K}(z,w).
$$
The canonical Hermitian structure of $E$ is then given by $H(z)=K(z,z)^{-1}$.

\subsection{\sf Induced representations} We briefly recall some known facts of representation theory.  
Let $G,H$ be real (or, complex) Lie
groups and $H\subseteq G$ be closed.  Given a representation $\varrho$ of $H$ on
a complex finite dimensional vector space $V$, let $\mathcal F(G, V)^H$ denote
the linear space of $C^\infty$ (or holomorphic) functions $F: G \to V$ satisfying
$$
F(g h) = \varrho(h)^{-1} F(g),\,\:\: g\in G,\, h \in H.
$$
The induced representation (cf. \cite[p. 187]{Kir}) $\mathbb U:=\mbox{\rm
Ind}_H^G(\varrho)$ acts on the linear space $\mathcal F(G, V)^H$ by left
translation:  $(\mathbb U_{g_1}f)(g_2) =  f(g_1^{-1}g_2)$.

From the linear representations $(\varrho, V)$ of $H$, one obtains all the $G$ -
homogeneous vector bundles  over $M=G/H$  as $G\times_H V$, which is $(G\times
V) / \sim$, where
$$
(gh, v) \sim (g, \varrho(h) v),\,\,h,g \in G,\,v\in V.
$$
The map $(g,v) \mapsto g H$ is clearly constant on the equivalence class
$[(g,v)]$ and hence defines a map $\pi: G\times_H V \to M$.  An action
$\hat{g}$, $g\in G,$ of the group $G$ is now defined on $G\times_H V$ by setting
$\widehat{g^\prime}\big ([(g,v)] \big ) = [(g^\prime g, v)]$.
This definition is independent of the choice of the representatives chosen.
Thus  $G\times_H V$ is a homogeneous vector bundle on $M$.
As in Section 1.1, there is a representation $U$ of $G$ on the
sections of $G\times_H V$, where $(U(g) s)(x) = \hat{g}
\big (s(g^{-1} \cdot x) \big )$. The lift to $G$ of the section $s$ of the
vector bundle $G\times_H V$ is $\tilde{s}: G\to V$ with $\tilde{s}(g) :=
\hat{g}^{-1} s(gH)$.  These again form the space $\mathcal F(G, V)^H$ which
shows that $U$ is just another realization of the representation $\mathbb U$.

When $M$ is a manifold with a group $G$ acting on it transitively, we use the 
usual identification $M=G/H$, where $H$ is the stabilizer in $G$ of a chosen fixed 
point $0 \in M$.  The map $q: g \mapsto g \cdot 0$ is the quotient map from $G$  to $M$.  
Suppose that there exists a global cross-section $p:M \to G$, that is, a map with $q \circ p = {\rm id}_{|\,M}$.  
Then $p$ gives a trivialization of the bundle $E=G\times_H V$.  
The trivializing map $\varphi$ is given for $v\in E_z$ by $\varphi(v)=(z,p(z)^{-1}v)$, that is, 
$\varphi_z=p(z)^{-1}$. (This $\varphi$ actually maps to $M \times E_0$, but $E_0$ with $H$ acting on 
it by the bundle action can be  identified with $(\varrho,V)$.)  As in Section 1.4, the action of 
$G$ on $E_z$ becomes $J_g(z)^{-1} = \varphi_{g z} \circ g_z \circ \varphi_z^{-1}$ which is now the 
group product $p(g z)^{-1}g p(z)$ (preserving the fibre $E_0$) followed by the identification of 
$E_0$ with $V=\mathbb C^n$, that is,
\begin{equation} \label{multformula}
J_g(z) = \varrho \big (  p(z)^{-1} g^{-1} p(g(z))\big ),\,\:\:z\in M,\, g\in G.
\end{equation}
The representation $U$ appears now as the multiplier representation with multiplier \eqref{multformula}.

Even though not needed in this paper, we point out that given any  $J:G\times M \to \mathrm{GL}_n(\mathbb C)$
satisfying the cocycle condition \eqref{multidnew}, the map $\big (U_{g}f \big ) (z) = J_g^{-1}(z)f(g^{-1} \cdot z)$ 
defines a multiplier representation of the group $G$.  Also,
it defines a representation $\varrho:h \mapsto J_{h^{-1}}(0)$ of the group $H$
on the vector space $V$.  The representation induced by this $\varrho$ is equivalent to $U$. In fact, 
the multiplier corresponding to the cross section $p$ and the
representation $\varrho$ is
\begin{eqnarray*}
\varrho\big ( p(z)^{-1}g^{-1}p(g\cdot  z)\big ) &=& J_{p(g\cdot z)^{-1} g
p(z)}(0)\\
&=& J_{p(z)}(0) J_{p(g\cdot z)^{-1} g}(p(z) \cdot 0)\\
&=& J_{p(z)}(0)J_g(z) J_{p(g\cdot z)^{-1}}(g \cdot z)\\
&=& J_{p(z)}(0)J_g(z)J_{p(g\cdot 0)}(0)^{-1}.
\end{eqnarray*}
which is equivalent to the multiplier $J$.

We remark that the inducing construction always gives a multiplier such that
$J_g(z) \in \varrho(H)$ for all $g, z$.  Not all multipliers possess this
additional property.  However, given any multiplier $J$, we can always find
another multiplier $J^\prime$ equivalent to $J$ such that $J_g^\prime(z) \in
\varrho(H)$, where $\varrho(h)=J_{h^{-1}}(0)$.  This is achieved by taking any
section $p$ and setting
$$
J_g^\prime(z) = J_{p(z)}(0) J_g(z) J_{p(g\cdot z)}(0)^{-1}.
$$

Holomorphic induced representation is a refinement of the induced representation
in the case of real groups $G, H$  such that $G/H$ has a $G$ - invariant complex
structure. The complex structure determines a subalgebra $\mathfrak b$ of
$\mathfrak g^\mathbb C$, namely the isotropy algebra in the local action of
$\mathfrak g^\mathbb C$ on $G/H$.  The holomorphic induced representation is the
restriction of the induced representation to a subspace of $\mathcal F(G,V)^H$,
defined by the differential equations $XF =  - \varrho(X)F$ for all 
$X\in\mathfrak b$, where $\varrho$ now is a representation of the pair $(H,\mathfrak b)$. 
It is an important fact that every $G$ - homogeneous holomorphic vector bundle arises by 
holomorphic induction from a simultaneous finite dimensional representation $\varrho$ of $H$ and $\mathfrak b$ 
(cf. \cite[Ch. 13]{Kir}). We will use this fact to determine all the holomorphic vector bundles which are 
homogeneous under the universal cover of the M\"{o}bius group.

\section{Homogeneous holomorphic vector bundles}
In the following, we explicitly construct all the irreducible homogeneous
holomorphic Hermitian vector bundles over the unit disc $\mathbb D$.  Every
homogeneous holomorphic Hermitian vector bundle on $\mathbb D$ is then obtained
as a direct sum of the irreducible ones (Corollary \ref{subhom}).  In Section 4,
we determine
which ones of these irreducible homogeneous holomorphic Hermitian vector bundles
 over $\mathbb D$ correspond to operators in the Cowen-Douglas class ${\mathrm
B}_n(\mathbb D)$.
\subsection{\sf The M\"{o}bius group}
Let $G_0$ be the M\"{o}bius group -- the group of bi-holomorphic automorphisms of 
the unit disc $\mathbb D$, $G = {\rm SU}(1,1)$ and $\mathbb K \subseteq G$ be the 
rotation group. Let $\tilde{G}$ be the universal covering group of $G$ 
(and also that of the group $G_0$).   
The group $G$ acts on the unit disc $\mathbb D$ according to the rule
$$g: z \mapsto (a z + b)(\bar{b} z + \bar{a})^{-1},\, g = 
\Big (\begin{matrix} a & b\\ \bar{b}& \bar{a} \end{matrix} \Big )\in G, \,\:\: z\in \mathbb D.$$
The group $\tilde{G}$ also acts on $\mathbb D$ (by $g\cdot z = q(g)\cdot z$, 
where $q:\tilde{G} \to G$ is the covering map), we denote the stabilizer of 
$0$ in it by $\tilde{\mathbb K}$.  So $\mathbb D \,\cong\, G/\mathbb K\, \cong\, \tilde{G}/\tilde{\mathbb K}$.
The complexification $G^{\mathbb C}$ of the group $G$ is the (simply connected) group ${\rm SL}(2,\mathbb C)$.

In the following, we use the notation of \cite{KM, KMISI}, which is the notation used in 
\cite{Sig}. The Lie algebra $\mathfrak g$ of the group $G$ is spanned by 
$X_1= \frac{1}{2} \mat{0}{1}{1}{0}$, $X_0 = \frac{1}{2}
\mat{i}{0}{0}{-i}$ and $Y= \frac{1}{2} \mat{0}{1}{1}{0}$.  The subalgebra
$\mathfrak k$ corresponding to $\mathbb K$ is spanned by $X_0$. In the
complexified Lie algebra $\mathfrak g^{\mathbb C}$, we  mostly use the complex
basis $h,x,y$ given by
\begin{eqnarray*}
h &=& - i X_0 = \frac{1}{2}\mat{1}{0}{0}{-1}\\
x &=& X_1+i Y = \mat{0}{1}{0}{0}\\
y &=& X_1 - i Y = \mat{0}{0}{1}{0}
\end{eqnarray*}


The subgroup of $G^{\mathbb C}$ corresponding to $\mathfrak g$ is $G$.
The group $G^{\mathbb C}$ has the closed subgroups
$\mathbb K^\mathbb C = \Big \{ \Big (\begin{matrix} z & 0\\ 0 & \tfrac{1}{z}
\end{matrix} \Big ):z \in \C, z\not = 0 \Big \}$, $P^+=\Big \{ \Big
(\begin{matrix} 1 & z\\ 0 & 1 \end{matrix} \Big ):z\in \C \Big \}$,
$P^- = \Big \{ \Big (\begin{matrix} 1 & 0\\ z & 1 \end{matrix} \Big ):
z\in \C \Big \}$; the corresponding Lie algebras
$\mathfrak k^\mathbb C = \Big \{ \Big (\begin{matrix} c & 0\\ 0 & -c
\end{matrix} \Big ): c\in \C \Big \}$, $\mathfrak p^+ = \Big \{ \Big
(\begin{matrix} 0 & c\\ 0 & 0 \end{matrix} \Big ):c\in \C \Big \}$,
$\mathfrak p^- = \Big \{ \Big (\begin{matrix} 0 & 0\\ c & 0
\end{matrix} \Big ):c\in \C \Big \}$
are spanned
by $h$, $x$ and $y$, respectively.  The product $\mathbb K^\mathbb C P^- =
\Big \{\mat{a}{0}{b}{\tfrac{1}{a}} : 0\not = a\in \C, b \in \C \Big
\}$ is a closed subgroup to be also denoted $B$; its Lie algebra is $\mathfrak b =
\mathbb C h + \mathbb C y$.  The product set $P^+\mathbb K^\mathbb CP^-=P^+B$
is dense open in $G^\mathbb C$, contains $G$, and the product decomposition of
each of its elements is unique. ($G^{\mathbb C}/B$ is the Riemann sphere, $g
{\mathbb K} \to g B,\:(g\in G)$ is the natural embedding of $\mathbb D\cong G/\mathbb K$ into it.)
Linear representations $(\varrho, V)$ of the algebra $\mathfrak b \subseteq
\mathfrak g^\mathbb C=\mathfrak s\mathfrak l(2,\mathbb C)$, that is, pairs of
linear transformations $\varrho(h), \varrho(y)$ satisfying
\begin{equation} \label{triangalgcommrel}
[\varrho(h),\varrho(y)] = - \varrho(y)
\end{equation}
are automatically representations of $\mathbb K$ as well.  Therefore they
give, by holomorphic induction, all the homogeneous holomorphic vector bundles.

A homogeneous holomorphic vector bundle that admits a $\tilde{G}$ - invariant
Hermitian structure will be called \emph{Hermitizable}.  Since the vector bundles 
corresponding to operators in the Cowen - Douglas class are of this type, we only 
consider these bundles (except for some comments following Remark \ref{Kob}). 
The $\tilde{G}$ - invariant Hermitian structures
on the homogeneous holomorphic vector bundle (making it into a homogeneous
holomorphic Hermitian vector bundle), if they exist, are given by
$\varrho(\tilde{\mathbb K})$ - invariant inner products on the representation
space $V$.  A $\varrho(\tilde{\mathbb K})$ - invariant inner product exists if
and only if $ \varrho(h)$ is diagonal with real diagonal elements in an
appropriate basis.   So, we will assume without restricting generality, that the
representation space of $\varrho$ is $\mathbb C^d$ and that $\varrho(h)$ is a
real diagonal matrix.

Furthermore, we will be interested only in irreducible homogeneous holomorphic
Hermitian vector bundles, this corresponds to $\varrho$ not being the orthogonal
direct sum of non-trivial representations.

Let $V_\lambda$ be the eigenspace of $\varrho(h)$ with eigenvalue $\lambda$.
We say that a Hermitizable homogeneous holomorphic vector bundle is
\emph{elementary} if the eigenvalues of $\varrho(h)$ form an uninterrupted string:
$-\eta, -(\eta+1),\ldots ,-(\eta+m)$.
Every irreducible homogeneous holomorphic Hermitian vector bundle is elementary,
In fact, let $-\eta$ be the largest eigenvalue of $\varrho(h)$ and $m$ be the
largest integer such that $-\eta, -(\eta+1), \ldots , -(\eta + m)$ are all
eigenvalues.  From \eqref{triangalgcommrel} we have $\varrho(y) V_\lambda
\subseteq V_{\lambda -1}$; this and orthogonality of the eigenspaces imply that
$V=\oplus_{j=0}^mV_{-(\eta +j)}$ and its orthocomplement are invariant under
$\varrho$.  So, $V$ is the whole space $\mathbb C^d$, and we have proved that the the bundle
is elementary.
We can write $V_{(\eta+j)}= \mathbb C^{d_j}$ and hence $(\varrho, \mathbb C^d)$ is described by the two matrices:
$$ \varrho(h) = \begin{pmatrix}
-\eta I_0 & &\\
&\ddots & \\
&&-(\eta + m) I_m
\end{pmatrix},$$
where $I_j$ is the identity matrix on $\mathbb C^{d_j}$ and
$$ Y:= \varrho(y) = \begin{pmatrix}
0 & & & \\
Y_1 & 0 & &\\
&Y_2&0 & &\\
&&\ddots&\ddots&\\
&&&Y_m&0
\end{pmatrix}$$
for some choice of matrices $Y_1, \ldots , Y_m$ that represent the linear transformations $Y_j:\mathbb C^{d_{j-1}} \to \mathbb
C^{d_j}.$ Let $E^{(\eta, Y)}$ denote the holomorphic bundle induced by this representation.

It is clear that $\varrho$ can be written as  the
tensor product of the one dimensional \rep $\sigma$ given by $\sigma(h) =
-\eta$, $\sigma(y)=0$, and the representation $\varrho^0$ given by
$\varrho^0(h)= \varrho(h) + \eta I$,
$\varrho^0(y) =\varrho(y)$.  Correspondingly, the bundle $E^{(\eta, Y)}$ for $\varrho$ is the
tensor product
of a line bundle $L_\eta$ and the bundle corresponding to $\varrho^0$, that is, $E^{(\eta, Y)} = L_\eta\otimes E^{(0,Y)}$.

For $g\in \tilde{G}$, $g^\prime(z)$ (we write $g^\prime(z) = \frac{\partial
g}{\partial z}(z)$) is a real analytic function on the
simply connected set $\tilde{G}\times \D$, holomorphic in $z$.
Also $g^\prime(z) \not = 0$ since $g$ is one-one and holomorphic.
Given any $\lambda \in \R$, taking the principal branch of the power function
when $g$ is near the identity, we can uniquely define $g^\prime(z)^\lambda$
as a real analytic function on $\tilde{G} \times \D$ which is holomorphic on
$\D$ for all fixed $g\in \tilde{G}$.

For the line bundle $L_\eta$, the multiplier is $g^\prime(z)^\eta$.
Consequently, the multiplier corresponding to the original $\varrho$ is
\begin{equation} \label{genmult}
J_g(z) = \big (g^\prime(z)\big )^\eta J_g^0(z),
\end{equation}
where $J^0$ is the multiplier obtained from $\varrho^0$.

The advantage of $\varrho^0$ is that it is also a representation of $G$ (not only of $\tilde{G}$) 
and extends to a representation of $G^\mathbb C$. The (ordinary) induced representation 
(in the holomorphic category) ${\rm Ind}_T^G(\varrho)$ operates on functions $F:G^\mathbb C \to V$ 
such that $F(g t)  = \varrho^0(t)^{-1}F(g)$ ($g\in G^\mathbb C$, $t\in T$).  
The restrictions of these functions $F$ to $G$ then give exactly the functions $\Phi:G \to V$ 
which satisfy $\Phi(g k) = \varrho^0(k)^{-1}\Phi(g)$ ($g\in G^\mathbb C$, $t\in T$) 
and $(X\Phi)(g) = - \varrho^0(X)\Phi(g)$ ($g\in G$, $X\in \mathfrak b$), that is, 
the space of the representation holomorphically induced by $\varrho^0$.  
Taking a holomorphic local cross section $p$ of $G^\mathbb C$ defined on $\mathbb D$, 
the functions $f(z)=F(p(z))$ give a trivialization of $E^{(0,Y)}$.

We use the local cross section $p:\D \to G^\bb{C}$, $z\mapsto p(z)
{:=}\Big (\begin{matrix} 1&z\\ 0&1
\end{matrix}\Big )$.
Apply \eqref{multformula} to compute the corresponding multiplier $J_g^0(z)$.
For $g=\Big (\begin{matrix} a & b\\ c & d
\end{matrix} \Big ) \in G$, we have 
\begin{eqnarray} \label{expY}
J_g^0(z)
&=& \varrho^0 \Big (
\Big (\begin{matrix} 1 & -g \cdot z\\ 0 & 1 \end{matrix} \Big )
\Big (\begin{matrix} a & b\\ c & d \end{matrix} \Big )
\Big (\begin{matrix} 1 & z\\ 0 & 1 \end{matrix} \Big ) \Big )^{-1} \nonumber\\
&=&\varrho^0  \begin{pmatrix}
c z+ d & 0 \\  -c & (c z+ d)^{-1}
\end{pmatrix} \nonumber\\
&=&\varrho^0 \Big (\,\Big (\begin{matrix} ( c z + d)^{\frac{1}{2}} & 0\\ 0 & (c
z + d)^{-\frac{1}{2}} \end{matrix} \Big ) \Big (\begin{matrix} 1 & 0\\ -c & 1
\end{matrix} \Big )
\Big (\begin{matrix} (c z + d)^{\frac{1}{2}} & 0\\ 0 & (c z + d)^{-\frac{1}{2}}
\end{matrix} \Big )\,\Big )\nonumber \\
&=& \varrho^0 (\exp(2 \log  (c z + d )^{\frac{1}{2}}h ))\varrho^0(\exp(-c
y))\varrho^0 (\exp(2 \log (c z + d )^{\frac{1}{2}} h) )\nonumber\\
&=&D_g(z) \exp(-c Y) D_g(z),
\end{eqnarray}
where $D_g(z)$ is the diagonal matrix with
$$
D_g(z)_{\ell \ell}=(c z +  d)^{-\frac{j}{2}}I_{d_j}.
$$

Computing the matrix entries of the exponential using \eqref{genmult}, we obtain for $g\in \tilde{G}$, $z\in \mathbb D$,
\begin{eqnarray} \label{multcovering}
(\!\!( J^{(\eta,Y)}_g(z))\!\!)_{p,\ell} &:=& (g^\prime(z))^\eta
J^0_g(z)\nonumber\\
&=&\begin{cases}
\frac{1}{(p - \ell)!} (-c_g)^{p-\ell}
g^\prime(z)^{\eta + \frac{p + \ell}{2}} Y_p \cdots Y_{\ell+1} &{~if~} p \geq
\ell\\
0&{~if~} p < \ell. \end{cases}
\end{eqnarray}
In this formula $c_g$ for $g \in \tilde{G}$  is to be understood as $c_{g^\#}$, where $g^\#$ is the 
projection to $G$ of $g$.
We note here, for later use, that there is also another way to interpret $c_g$ for $g\in \tilde{G}$. 
Taking a small neighborhood $\tilde{U}$ of the identity in $\tilde{G}$ such that the projection is a 
diffeomorphism onto a neighborhood $U$ of the identity in $G$, by computing in $U$, we find that
\begin{equation}\label{doubleder}
g^{\prime \prime}(z)=-2c_g\, g^\prime(z)^{3/2}
\end{equation}
holds with $c_g$ an analytic function of $g$ on $U$, independent of $z$.  
This is then true for $g\in \tilde{U}$ and by analytic continuation for all 
$g\in \tilde{G}$. So the equation \eqref{doubleder} can serve as a definition for $c_g$.


\begin{Proposition} \label{elem} All elementary Hermitizable homogeneous holomorphic
vector bundles are of the form $E^{(\eta,Y)}$ with $\eta \in \mathbb R$ and $Y$ as before.
The bundles $E^{(\eta, Y)}$ and $E^{(\eta^\prime, Y^\prime)}$ are isomorphic if and only if 
$\eta=\eta^\prime$ and $Y^\prime = AYA^{-1}$ with a block diagonal matrix $A$.
\end{Proposition}

\begin{proof}
The induced bundles are isomorphic if and only if the inducing representations $\varrho$, 
$\varrho^\prime$ are linearly equivalent, that is, $\varrho^\prime = A \varrho A^{-1}$ for some $A$.  
Since we have normalized the representations by fixing the matrix $\varrho(h)$, 
the equivalence must be given by an $A$ which commutes with $\varrho(h)$, that is, by a block diagonal $A$.
\end{proof}

Thus $E^{(\eta,\{Y\})}= L_\eta \otimes E^{(\{Y\})}$  parametrizes the
equivalence classes of elementary Hermitizable homogeneous holomorphic vector
bundles.  Here, we have let $\{Y\}$ denote the conjugacy class of $Y$ under
conjugation by a block diagonal matrix $A$.
\subsection{\sf Homogeneous holomorphic Hermitian vector bundles}

We proceed to discuss homogeneous holomorphic Hermitian vector bundles.  From here on we  
will always use the trivialization we just described.  We will not always make a careful 
distinction between a section of $E^{(\eta,Y)}$ and the functions from $\mathbb D$ to 
$\mathbb C^d$on which $G$ acts by the multiplier $J^{(\eta,Y)}_g(z)$. As in Section 1, 
a Hermitian structure appears in the trivialization as a family of quadratic forms 
$\inner{H(z)\xi}{\xi}$, which because of the homogeneity is determined by a single 
positive definite block-diagonal matrix $H=H(0)$.  We denote by $(E^{(\eta, Y)}, H)$ the bundle 
$E^{(\eta, Y)}$ equipped with the Hermitian structure determined by $H$.

\begin{Proposition} \label{prop2.2}
The Hermitian vector bundles $(E^{(\eta, Y)}, H)$ and $(E^{(\eta^\prime, Y^\prime,)}, H^\prime)$ are isomorphic if and only if
$\eta=\eta^\prime$, $Y^\prime = AYA^{-1}$ and $H^\prime = {A^*}^{-1} H A$ with a block diagonal matrix $A$.
\end{Proposition}
\begin{proof}
The trivialization of the sections obtained by starting with $\varrho$ (resp. $\varrho^\prime = A\varrho A^{-1}$) 
are related as $f^\prime (z) =  Af(z)$.  Now,  $H^\prime(z)$ gives the same metric as $H(z)$ if and only if 
$\inner{H^\prime(z) f^\prime(z)}{ f^\prime(z)} = \inner{H(z)f(z)}{f(z)}$.   From this, the statement follows.
\end{proof}

For any $H$, clearly there is an $A$ such that ${A^*}^{-1}HA = I$. This means that every elementary homogeneous 
holomorphic Hermitian vector bundle is isomorphic to one of the form $(E^{(\eta, Y}, I)$.  Two vector bundles of 
this form are equivalent if and only if $Y^\prime = A Y A^{-1}$ with $A$ such that ${A^*}^{-1} I A^{-1} = I$, 
that is, with a block-diagonal unitary $A$.  We denote by $[Y]$ the equivalence class of $Y$ under conjugation 
by block-diagonal unitaries and write $E^{(\eta, [Y])}$ for the equivalence class of $(E^{(\eta, Y)}, I)$, 
omitting the $I$.  We now have the first half of the following Proposition. 
\begin{Proposition}  \label{prop2.3}
For $\eta \in \mathbb R$, $[Y]$ a block unitary conjugacy class of matrices $Y$, the vector bundles $E^{(\eta, [Y])}$ 
form a parametrization of the elementary homogeneous holomorphic Hermitian vector bundles. The Hermitian vector bundle  
$E^{(\eta, [Y])}$ is irreducible if and only if $Y$ cannot be split into orthogonal direct sum 
$Y^\prime \oplus Y^{\prime\!\prime}$ with $Y^\prime$, $Y^{\prime\!\prime}$ having the same block diagonal form as $Y$.
\end{Proposition}
\begin{proof}
The last statement follows since the irreducibility of $E^{(\eta, [Y])}$ is the same as the possibility of 
splitting $\varrho$ into an orthogonal direct sum of two sub-representations.
\end{proof}
Proposition \ref{prop2.3}, with a different proof, also appears in \cite{BisGM}.

The following Theorem is important because its hypothesis is exactly what we know about the vector bundle 
corresponding to a homogeneous operator in the Cowen-Douglas class ${\rm B}_n(\mathbb D)$.  It was stated 
in \cite{KM} but proved without the uniqueness statement. Here we give a complete proof. 

\begin{Theorem} \label{uniqlift}
Let $E$ be a Hermitian holomorphic vector bundle over $\mathbb D$ and suppose
that for all $g \in G$, there exists an automorphism of $E$ whose action on
$\mathbb D$ coincides with $g$.  Then the full automorphism group of $E$ is   
reductive and $\tilde{G}$ acts on $E$ by automorphisms
in a unique way.
\end{Theorem}

\begin{proof}
Let $\hat{G}$ denote the connected component of the automorphism group of $E$.
It is a Lie group because it is the connected component of the isometry group of
$E$ under the Riemannian metric defined for vectors tangent to the fibres by the
Hermitian structure and for vectors horizontal with respect to the Hermitian
connection by the $G$-invariant metric of $\mathbb  D$.

Let $N \subseteq \hat{G}$ be the subgroup of elements acting on $\mathbb D$ as
the identity map.   The subgroup $N$ is normal, and the projection $\pi:\hat{G}
\to G$ is a homomorphism with kernel $N$.  Let $\mathbb K$ be the stabilizer of $0$ in
$G$ and let $\hat{\mathbb K} = \pi^{-1}(\mathbb K)$. The group $\hat{\mathbb K}$ contains $N$ and is
compact because it is the stabilizer of the origin in the fiber over $0$.

Let $\hat{\mathfrak g}, \mathfrak g, \mathfrak k, \mathfrak n, \hat{\mathfrak
k}$
be the Lie algebras of the groups defined above, and let $\mathfrak g =
\mathfrak k + \mathfrak p$ be the Cartan decomposition.  Since $\hat{\mathbb K}$ is
compact,  we can choose an $\mathrm{Ad}(\hat{\mathbb K})$ - invariant complement
$\hat{\mathfrak p}$ to $\hat{\mathfrak k}$ in $\hat{\mathfrak g}$.  Now, $\pi_*$
maps $\hat{\mathfrak k}$ onto $\mathfrak k$ with kernel $\mathfrak n$.  By
counting dimension,  it follows that $\pi_*$ maps $\hat{\mathfrak p}$ to
$\mathfrak p$ bijectively.

We set $\hat{\mathfrak k}_0 = [\hat{\mathfrak p}, \hat{\mathfrak p}]$.  Then
$\pi_*(\hat{\mathfrak k}_0) = [\pi_*\hat{\mathfrak p}, \pi_*\hat{\mathfrak p}] =
\mathfrak k$, therefore $\hat{\mathfrak k}_0 \subseteq \pi_*^{-1}(\mathfrak k) =
\hat{\mathfrak k}$.  It follows that $[\hat{\mathfrak k}_0, \hat{\mathfrak
p}]\subseteq \hat{\mathfrak p}$ and by the Jacobi identity, $\hat{\mathfrak
g}_0= \hat{\mathfrak k}_0 + \hat{\mathfrak p}$ is a subalgebra.  Similarly,
$[\mathfrak n, \hat{\mathfrak p}] \subseteq \hat{\mathfrak p}$ since
$\mathfrak n \subseteq \hat{\mathfrak k}$.  But $\mathfrak n$ is an ideal, so
$[\mathfrak n, \mathfrak p]=0$, and by the Jacobi identity $[\mathfrak n,
\hat{\mathfrak g}_0]=0$.  Finally, $\hat{\mathfrak g} = \mathfrak n \oplus
\hat{\mathfrak g}_0$ and $\mathfrak g$ is reductive.
The analytic subgroup $\hat{G}_0 \subseteq \hat{G}$  corresponding to
$\hat{\mathfrak g}_0$ is a covering group of $G$ and  therefore
it acts on $E$ by automorphisms.  It is the unique subgroup of $\hat{G}$ with
this property because $\hat{\mathfrak g}_0$, being the maximal semisimple ideal
in the reductive algebra $\hat{\mathfrak g}$,  is uniquely determined. The
action of $\hat{G}_0$ now lifts to a unique action of $\tilde{G}$.
\end{proof}
Theorem \ref{uniqlift} implies that every homogeneous operator in
the Cowen-Douglas class $\mathrm{B}_n(\mathbb D)$ has an associated
representation irrespective of whether it is irreducible or not.  The following
Corollary has also appeared in \cite{BisGM}.
\begin{Corollary} \label{subhom}
If a  Hermitian holomorphic vector bundle $E$ is homogeneous  and is reducible
($E=E_1\oplus E_2$)  as a Hermitian holomorphic vector bundle then it is
reducible as a homogeneous Hermitian holomorphic vector bundle, that is, $E_1$
and $E_2$ are also homogeneous.
\end{Corollary}
\begin{proof} We consider the automorphisms $\exp th$ of $E$, where
$$h = \begin{cases}i I & \mbox{ on } E_1 \\
-i I & \mbox{ on } E_2.
    \end{cases}$$
Then $h$ is in $\mathfrak n$ since $\exp th$ ($t\in \mathbb R$) preserves fibres.  So, $h$
commutes with $\hat{\mathfrak g}_0$.  The sections of $E_1, E_2$ are
characterized as eigensections of $h$ corresponding to different eigenvalues. Thus $\hat{G}_0$, and its
universal covering $\tilde{G}$  preserve the eigensections of $h$.
\end{proof}

\section{Homogeneous holomorphic Hermitian vector bundles with reproducing kernel}

In this Section, we determine which ones of the elementary homogeneous holomorphic
Hermitian vector bundles have their Hermitian structure coming from a reproducing kernel.
In other words, which are the homogeneous holomorphic vector bundles that have a $\tilde{G}$ - 
invariant reproducing kernel $K(z,w)$.
When there is a reproducing kernel $K$, it
gives a canonical Hermitian structure by setting $H=K(0,0)^{-1}$.
Let $p_z =\frac{1}{\sqrt{1-|z|^2}} \left (\begin{smallmatrix} 1 & z
\\ \bar{z} & 1 \end{smallmatrix} \right ) \in  G$, so $ p_z \cdot 0 = z.$
Writing $J_z^{(\eta, Y)}$ for $J_{p_z}^{(\eta,Y)}(z)$, we have
\begin{equation} \label{tranrulK}
K(z,z) = J_z^{(\eta,Y)} K(0,0) {J_z^{(\eta,Y)}}^*.
\end{equation}
So, the question amounts to enumeration of all the possibilities for $K(0,0)$.

\subsection{\sf An intertwining map}
For $\lambda >0$, let $\mathbb A^{(\lambda)}$ be the Hilbert space of holomorphic
functions on the unit disc with reproducing kernel $(1-z \bar{w}))^{-2 \lambda}$.
It corresponds to the homogeneous line bundle $L_\lambda$.  The group $\tilde{G}$ 
acts on it unitarily with the multiplier $g^\prime(z)$.  This action is the 
Discrete series representation $D^{(\lambda)}_g$. 
Let $\mathbb C^d = \oplus_{j=0}^{d_j} \mathbb C^{d_j}$. We think of functions 
$f: \mathbb D \to \mathbb C^d$ as having components $f_j: \mathbb D \to \mathbb C^{d_j}$.  
Let $\mathbf A^{(\eta)}= \oplus_{j=0}^m \mathbb A^{(\eta +
j)} \otimes \mathbb C^{d_j}$.  For $\eta > 0$, $Y$ as before and 
$f_j \in \mathbb A^{(\eta + j)} \otimes \mathbb C^{d_j}$, define
\begin{equation} \label{Gammadefn}
\big (\Gamma^{(\eta, Y)} f_j \big )_\ell =\begin{cases}
\frac{1}{(\ell -j)!}\frac{1}{(2 \eta + 2 j)_{\ell -j}} Y_\ell \cdots Y_{j+1}
 f_j^{(\ell -j)}&{~if~} \ell \geq j\\
0&{~if~}  \ell < j.\end{cases}
\end{equation}
So, $\Gamma^{(\eta, Y)}$ maps ${\rm Hol}(\mathbb D, \mathbb C^d)$ into itself.
Let $N$ be an invertible $d\times d$ block diagonal matrix with blocks $N_j$, $0\leq j
\leq m$, $d=d_0+\cdots +d_m$.  We will assume throughout that $N_0=I_{d_0}$.
This is only a normalizing condition.
We  can normalize further by assuming that each block diagonal matrix with 
$d_j \times  d_j$ blocks $N_j$ is positive definite but this is not important.
We can think of $N$ as changing the natural inner product of each 
$\mathbb C^{d_j}$ to $\inner{N_ju}{N_jv}_{\mathbb C^{d_j}}$.
We  let $\Gamma_N^{(\eta, Y)} = \Gamma^{(\eta, Y)}\circ N$ and
$\mathcal H_N^{(\eta, Y)}$ denote the image of  $\Gamma_N^{(\eta, Y)}$ in the
space of holomorphic functions  $\mathrm{Hol}(\mathbb D, \mathbb C^d)$.

\begin{Theorem} \label{equiv}
The map $\Gamma_N^{(\eta, Y)}$ is a $\tilde{G}$ - equivariant isomorphism of $\mathbf A^{(\eta)}$ onto the Hilbert space 
$\mathcal H_N^{(\eta, Y)}$ on which the $\tilde{G}$
action is unitary via the multiplier $J_g^{(\eta, Y)}(z)$. It has a reproducing kernel $K^{(\eta, Y)}_N(z,w)$ such that
$$
\big (K^{(\eta, Y)}_N(0,0) \big )_{\ell \ell}=\sum_{j=0}^\ell \frac{1}{(\ell
-j)!}\frac{1}{(2 \eta + 2j)_{\ell -j}} Y_\ell \cdots Y_{j+1} N_jN_j^*
Y_{j+1}^*\cdots Y_\ell^*.
$$
\end{Theorem}
\begin{proof}
The injectvity of the map $\Gamma_N^{(\eta, Y)}$ is clear from its definition.  
It is also apparent that the image $\mathcal H_N^{(\eta, Y)}$ is the algebraic direct sum of the summands
$\mathbb A^{(\eta+j)}\otimes \mathbb C^{d_j}$ of $\mathbf A^{(\eta)}$.  We define a norm on
$\mathcal H_N^{(\eta, Y)}$ by stipulating that $\Gamma_N^{(\eta, Y)}$ is a Hilbert space isometry.
This gives us the Hilbert space  $\mathcal H_N^{(\eta, Y)}$ and the unitary action $U_g$ of $\tilde{G}$ on it.  
We have to show that it is the multiplier action given by $J_g^{(\eta, Y)}(z)$. For this, we must verify that


\begin{equation} \label{verify}
\Gamma_N^{(\eta,Y)} \circ\big (\oplus d_j D_{g^{-1}}^{(\eta+j)} \big )
=U_{g^{-1}}\circ \Gamma_N^{(\eta, Y)}.
\end{equation}
Since $N$ obviously intertwines $\oplus\,d_jD^{(\eta+j)}$ with itself, it
suffices to prove \eqref{verify} for $\Gamma^{(\eta,Y)}$ in place of
$\Gamma^{(\eta,Y)}_N=\Gamma^{(\eta,Y)}\circ N$.  Furthermore,
it is enough to verify this relation for each $f\in \mathbb A^{(\eta+j)} \otimes
\mathbb C^{d_j}$, that is, to show
$$
\Gamma^{(\eta, Y)}\big ((g^\prime)^{\eta+j}(f\circ g)\big ) = J_g\big
((\Gamma^{(\eta, Y)} f) \circ g\big ), \, f \in \mathbb A^{(\eta+j)} \otimes
\mathbb C^{d_j},\, 0\leq j \leq m.
$$
We will show that the $\ell$th components on both sides are equal.
First, if $\ell < j$ then both sides are $0$.  Second if $\ell \geq j$, on the
one hand, using Lemma 3.1 of \cite{KM}  which  is easily proved by induction starting from
the equation \eqref{doubleder} and says that
\begin{equation}\label{Leibnitz}
\big ( (g^\prime)^\ell(f\circ g) \big )^{(k)} =
\sum_{i=0}^k
\tbinom{k}{i}(2\ell+i)_{k-i}(-c)^{k-i}(g^\prime)^{\ell+\frac{k+i}{2}}\big (
f^{(i)} \circ g \big )
\end{equation}
for any $g\in \tilde{G}$, we have
\begin{eqnarray*}
\lefteqn{
\Gamma^{(\eta, Y)}\big ((g^\prime)^{\eta+j}(f\circ g)\big )}\\
&=& \frac{1}{(\ell -j)!}\frac{1}{(2 \eta + 2 j)_{\ell -j}} Y_\ell \cdots Y_{j+1}
\big ( (g^\prime)^{\eta+j} (f\circ g) \big )^{(\ell-j)}\\
&=& \frac{1}{(\ell -j)!}\frac{1}{(2 \eta + 2 j)_{\ell -j}} Y_\ell \cdots Y_{j+1}
\sum_{i=0}^{\ell-j} \tbinom{\ell -j}{i}
(2\eta+2j+i)_{\ell-j-i}(-c)^{\ell-j-i}(g^\prime)^{\eta+j+\frac{\ell-j+i}{2}}
(f^{(i)})\circ g\\
&=&  Y_\ell \cdots Y_{j+1} \sum_{i=0}^{\ell-j} \frac{1}{(\ell -j-i)!i!}
\frac{1}{(2\eta+2j)_i}(-c)^{\ell-j-i}(g^\prime)^{\eta+j+\frac{\ell-j+i}{2}}
(f^{(i)})\circ g,
\end{eqnarray*}
On the other hand,
\begin{eqnarray*}
\lefteqn{\sum_{p=j}^m (J_g)_{\ell p} \big ((\Gamma^{(\eta, Y)} f)_p \circ g\big
)}\\
&=& \sum_{p=j}^\ell
(-c)^{\ell-p}\frac{1}{(\ell-p)!}(g^\prime)^{\eta+\frac{p+\ell}{2}}Y_\ell \cdots
Y_{p+1}\frac{1}{(p -j)!}\frac{1}{(2 \eta + 2 j)_{p -j}} Y_p \cdots Y_{j+1}
f^{(p-j)}\circ g\\
&=& \sum_{p=j}^\ell \frac{1}{(\ell-p)!}\frac{1}{(p -j)!}\frac{1}{(2 \eta + 2
j)_{p -j}}(-c)^{\ell-p}(g^\prime)^{\eta+\frac{p+\ell}{2}}Y_\ell \cdots Y_{j+1}
f^{(p-j)}\circ g\\
&=& \sum_{i=0}^{\ell-j} \frac{1}{(\ell-i-j)!}\frac{1}{i!}\frac{1}{(2 \eta + 2
j)_{i}}(-c)^{\ell-j-i}(g^\prime)^{\eta+\frac{j+i+\ell}{2}}Y_\ell \cdots Y_{j+1}
f^{(p-j)}\circ g.
\end{eqnarray*}
This completes the verification of \eqref{verify}.
Finally, we we observe that  $\mathcal H_N^{(\eta, Y)}$ has a reproducing kernel 
$K^{(\eta, Y)}_N(z,w)$ because it is  the image of $\mathbf A^{(\eta)}$ under an 
isomorphism given by a holomorphic differential operator, so point evaluations 
remain continuous.  Then  $K^{(\eta, Y)}_N(z,w)$ is obtained by applying 
$\Gamma_N^{(\eta, Y)}$ to the reproducing kernel of $\mathbf A^{(\eta)}$ 
once as a function of $z$ and once as a function of $w$. This computation 
is easily carried out and gives the formula for $K^{(\eta, Y)}_N(0,0)$.
\end{proof}
Writing $H:=H^{(\eta, Y)}_N = \big (K_N^{(\eta, Y)}(0,0)\big )^{-1}$, the Hilbert
space $\mathcal H_N^{(\eta, Y)}$ is a space of sections of the homogeneous
holomorphic Hermitian vector bundle $(E^{(\eta, Y)},H)$ in our
(canonical) trivialization.
\begin{Theorem}  \label{3.1}
The construction with $\Gamma_N^{(\eta, Y)}$ gives all elementary homogeneous 
holomorphic Hermitian vector bundles which possess a reproducing kernel, namely, those of the form
$$\big (E^{(\eta, Y)}, (K^{(\eta, Y)}_N(0,0))^{-1} \big ),$$
where $\eta > 0$, $Y$ are arbitrary and $K^{(\eta, Y)}_N(0,0)$ is the form given in 
Theorem \ref{equiv}.  The vector bundles $\big (E^{(\eta, Y)}, (K^{(\eta, Y)}_N(0,0))^{-1} \big )$
and $\big (E^{(\eta^\prime, Y^\prime)}, (K^{(\eta^\prime, Y^\prime)}_{N^\prime}(0,0))^{-1} \big )$ 
are equivalent if and only if $\eta = \eta^\prime$, 
$Y^\prime = AY A^{-1}$ and $N^\prime {N^\prime}^* = ANN^*A^*$ for some 
invertible block diagonal matrix $A$ of size $d\times d$.
\end{Theorem}

\begin{proof}
The existence of a reproducing kernel  implies that the vector bundle is
Hermitizable.
Such a bundle is of the form $(E^{(\eta, Y)}, H)$ by Propositions \ref{elem} and \ref{prop2.2}.
When it has a reproducing kernel, then in our canonical trivialization this is a
matrix valued function $K(z,w)$, and we have $H=K(0,0)^{-1}$. The $\tilde{G}$
action $U$ which is now unitary, is given by the multiplier $J_g^{(\eta, Y)}(z)$.
The equation \eqref{multcovering} shows that the action of $\tilde{\mathbb K}$ is diagonalized
by the polynomial vectors:  If $v_j \in \mathbb C^{d_j}$ and $f(z) = z^n v_j$, 
then for $k_\theta$ such that $k_\theta (z) = e^{i \theta} z$, we have $U_{k_\theta} f = e^{i \theta(\eta+j +k)}f$.
It follows that $U$ is a direct sum of the Discrete  series representations
$D^{(\eta+j)}$, $0\leq j \leq m$.  In particular, it follows that $\eta > 0$.

The map $\Gamma^{(\eta, Y)}$ (and $\Gamma_N^{(\eta, Y)}$ for any block diagonal
$N$) intertwines the representations $U$ and $ \oplus_{j=0}^m \,d_j D^{(\eta+j)}$,  
both of which are unitary.  By Schur's Lemma it follows that $N$ can be
chosen  such that $\Gamma^{(\eta, Y)}\circ N$ is unitary.  This proves that the
bundle $E^{(\eta, Y)}$ corresponds to the Hilbert space $\mathcal H_N^{(\eta, Y)}$.

The statement about equivalence follows from the analogous  statement in Proposition \ref{prop2.2}. 
\end{proof}

\begin{Remark}\label{Kob}{\em In the proof we only used the unitarizability of the
$\tilde{G}$ action
on the sections of the Hermitizable bundle $E$.  In this vein, an even more
general result holds:

For any $\tilde{G}$ - homogeneous holomorphic vector bundle $E$, if the
$\tilde{G}$ action  on the sections is unitarizable then it is a direct sum of
bundles corresponding to some  $\mathcal H_N^{(\eta, Y)}$.
(The possible unitary structures correspond to different choices of $N$.)}
\end{Remark}

One way to prove this is to use the ``Inverse propagation theorem'' of T.
Kobayashi \cite{TK}.   If the action of $\tilde{G}$ is unitary, then so is the
$\tilde{\mathbb K}$ action on the fibres,  and we are back
in the situation of Theorem \ref{3.1}.

Here we sketch a more direct proof which also shows what the non-Hermitizable 
homogeneous holomorphic vector bundles are like.

A general $E$ is still gotten from two matrices $Z=\varrho(h),\, Y= \varrho(y)$
such that $[Z,Y]=-Y$ by holomorphic induction.  The inclusion $YV_\lambda
\subseteq V_{\lambda -1}$ still holds for the generalized eigenspaces of $Z$.
Using some easy identities for $g^\prime(z)$, we can then verify that
$$
J_g(z) = \exp\big (\textstyle{\frac{1}{2}} (\log(g^\prime(z))^\prime Y\big
)\exp\big (-\log g^\prime(z) Z\big ),
$$
which, in the case where $Z$ is diagonal, is just another way to write \eqref{Gammadefn}, is a multiplier.

Writing $U_g$ for the action of $\tilde{G}$ on $\mathrm{Hol}(D,V)$ given by
$J_g(z)$,
we compute
\begin{equation} \label{J_g}
\big (U_{\exp t i h}f\big )(z)= \exp(i t Z)f(e^{-it}z).
\end{equation}
Hence $(U_hf)(0)= Zf(0)$ and by a similar computation $(U_yf)(0)= Yf(0)$.  This
shows that $J_g(z)$ gives a trivialization of our $E$.  It also shows that
$U_k,\, k\in \tilde{\mathbb K}$ maps the spaces $\mathcal M_p$ of monomials of
degree $p$ to $\mathcal M_p$ for all $p \geq 0$.  Hence $\tilde{\mathbb K}$ -
finite vectors are exactly the ($V$ - valued) polynomials.

Now if $U$ is unitary with respect to some inner product, then it is a sum of
irreducible representations. The $\tilde{\mathbb K}$-types of these (i.e. the
eigenfunctions of $U_h$) are known to be one dimensional  and together they
span
the space of $\tilde{\mathbb K}$ - finite vectors.  By  \eqref{J_g}, $U_h$ maps
any $z^p v \in\mathcal M_p$ to $z^p(Zv - pv)$.  It follows that $Z$ must be
diagonalizable, otherwise the eigenfunctions of $U_h$ could not span $\mathcal
M_p$.

\subsection{\sf Parametrization} The description of the homogeneous 
holomorphic Hermitian vector bundles given in Theorem \ref{3.1} can be 
made more explicit.  We now proceed to determine, in terms of the parametrization 
$E^{(\eta, [Y])}$ of elementary homogeneous holomorphic Hermitian vector bundles 
as in Proposition \ref{prop2.3}, exactly which ones of these have their Hermitian 
structure come from a reproducing kernel.
\begin{Theorem} \label{Pcond}
The Hermitian structure of $E^{(\eta, [Y])}$ comes from a ($\tilde{G}$ - invariant)
reproducing kernel if and only if $\eta > 0$ and
\begin{equation*}
 I - Y_j \big (\sum_{k=0}^{j-1}\frac{(-1)^{j+k}}{(j-k)!(2\eta + j+k - 1)_{j-k}} Y_{j-1}\cdots Y_{k+1} Y_{k+1}^*
\cdots Y_{j-1}^* \big )Y_j^* > 0
\end{equation*}
for $j=1,2,\ldots ,m$.
\end{Theorem}
\begin{proof}
We have a description of all the vector bundles with reproducing kernel in Theorem \ref{3.1}.  
To see how  this description appears in the parametrization $E^{(\eta, Y)}$, 
we have to answer the question:  For what $\eta, [Y]$, is it possible to find a block-diagonal $N$ 
such that $K_N^{(\eta, Y)}(0,0) = I$.  Writing this more explicitly, we have the system of equations
\begin{equation} \label{P}
I_\ell - \sum_{j=0}^{\ell} \frac{1}{(\ell
-j)!}\frac{1}{(2\eta + 2j)_{\ell-j}} Y_\ell\cdots Y_{j+1} N_j N_j^*Y_{j+1}^*
\cdots Y_{\ell}^* = 0,
\end{equation}
$\ell = 1,\ldots , m$ and the question is if the solution $N_jN_j^*$, $j=1,\ldots , m$ 
consists of positive definite matrices.

We claim that the solution of \eqref{P} is given by
\begin{equation} \label{Betafn}
 N_jN_j^* = \sum_{k=0}^{j}\frac{(-1)^{j+k}}{(j-k)!(2\eta + j+k - 1)_{j-k}} Y_j \cdots Y_{k+1} Y_{k+1}^*
\cdots Y_j^*,
\end{equation}
for $j=1,\ldots , m$.

In fact, substituting the expression for $ N_jN_j^*$ from \eqref{Betafn} into \eqref{P}, we have
$$
I_\ell - \sum_{j=0}^\ell \sum_{k=0}^j  \frac{1}{(\ell-j)!(2\eta+2j)_{\ell-j}}
\frac{(-1)^{j+k}}{(j-k)!(2\eta+j+k-1)_{j-k}} Y(k) = 0,
$$
where $Y(k) =  Y_\ell \cdots Y_{k+1} Y_{k+1}^*\cdots Y_\ell^*$.  The coefficient of $Y(k)$, from the above, is seen to be
$$
\frac{1}{(\ell-k)!^2}\sum_{j=k}^\ell (-1)^{j+k} \binom{\ell-k}{j-k} (2\eta+2j -1) B(2\eta+k+j-1, \ell-k+1),
$$
where $B(x,y) = \frac{\Gamma(x)\Gamma(y)}{\Gamma(x+y)}$ is the usual Beta function.
Using the binomial formula and the integral representation: $B(x,y) = \int_0^1 t^{x-1}(1-t)^{y-1}dt$, 
it simplifies to  
\begin{eqnarray*}
\lefteqn{\frac{1}{(\ell-k)!^2} \int_0^1 \Big\{ (2\eta+2k -1)t^{2\eta+2k-2}(1-t)^{2(\ell-k)} - 
2(\ell-k) t^{2\eta+2k-1}(1-t)^{2(\ell-k)-1}\Big \} dt}\\
&=&\frac{1}{(\ell-k)!^2} \int_0^1 \Big\{t^{2\eta+2k-2}(1-t)^{2(\ell-k) - 1} \big ( (2\eta + 2k -1)  
- \big ( 2\eta  + 2\ell -1) t\big )\Big \} dt\\
&=& \frac{1}{(\ell-k)!^2} \big (xB(x,y)-(x+y)B(x+1,y)\big ),
\end{eqnarray*}
where $x=2\eta+2k-1$ and $y=2\ell-2k$, which is zero except when $k = 0 = \ell$. This verifies the claim.

The right hand side of the equation \eqref{Betafn} is exactly the expression given in the statement 
of the Theorem, so its positivity is the condition we were seeking.
\end{proof}

When $Y$ is given, we may ask what are the values of $\eta$ for which the positivity 
condition of the Theorem holds. It obviously holds when $\eta$ is large. We can also see that
there exists a number $\eta_Y > 0$ such that it holds if and only if $\eta > \eta_Y$.
For this we only have to see that if $E^{(\eta, Y)}$ has a reproducing kernel for 
some $\eta > 0$, then so does $E^{(\eta + \varepsilon, Y)}$ for all $\varepsilon > 0$.  
Now $E^{(\eta + \varepsilon, Y)} = L_\varepsilon \otimes E^{(\eta, Y)}$ which shows that 
the product $K(z,w) = (1-z\bar{w})^{-2 \varepsilon} K_I^{(\eta, Y)}(z,w)$ is a reproducing 
kernel for $E^{(\eta + \varepsilon, Y)}$, and $K(0,0) =I$ still holds.


When $m=1$, the condition of the Theorem \ref{Pcond} reduces to
$$
I - \frac{1}{\eta} Y_1Y_1^* > 0.
$$
In this case, we have $\eta_Y = \frac{1}{2} \|Y_1Y_1^*\|$ in terms of the usual matrix norm.

\section{\sf Classification of the homogeneous operators in 
the Cowen-Douglas class}
The following theorem together with Theorems \ref{equiv} and \ref{3.1}, and Corollary 
\ref{subhom} gives a complete classification of homogeneous operators in the Cowen-Douglas class.

\begin{Theorem}
All the homogeneous holomorphic Hermitian vector bundles with a reproducing
kernel correspond to homogeneous operators in the Cowen-Douglas class. The
irreducible ones are the adjoint of the multiplication operator $M$ on the
space  $\mathcal H_I^{(\eta, Y)}$ for some  $\eta >0$ and irreducible  $Y$.
The block matrix { $Y$} is determined up to conjugacy by block diagonal
unitaries.
\end{Theorem}
\begin{proof}
First we note that by Theorems \ref{3.1} and \ref{Pcond} every homogeneous holomorphic  Hermitian vector bundle 
can be written in the form $(E^{\eta, Y},I)$ with $\eta >0$. The Hilbert space $\mathcal H^{\eta, Y}_I$ is a 
subspace of the  (trivialized) holomorphic sections of $(E^{\eta, Y},I)$ which is the image under the map 
$\Gamma^{(\eta,Y)}_N$ of $\mathbf A^{(\eta)}$.  We have to show only that the operator $M^*$ on $\mathcal H^{(\eta, Y)}_I$ 
belongs to the Cowen-Doulas class.  For this we compute the matrix of $M$ in an appropriate orthonormal basis.

Each of the Hilbert spaces $\mathbb A^{(\eta+j)}$ ($0\leq j \leq m$) has a natural orthonormal basis
$$\Big \{e_j^n(z):=\sqrt{\frac{(2\eta+2j)_n}{n!}} z^n\,:n \geq 0 \Big \}.$$

Hence $\mathbb A^{(\eta+j)} \otimes \mathbb C^{d_j}$ has the basis $e_j^n \varepsilon_{q}^{(j)}$, 
where $\{\varepsilon_q^{(j)}: 1 \leq q \leq d_j\}$ is the natural basis of $\mathbb C^{d_j}$.  The Hilbert space $\mathbf
A^{(\eta)}$ then has the orthonormal basis 
$e_j^n \varepsilon_{jq}$ with $\varepsilon_{jq} := \varepsilon_j \otimes \varepsilon^{(j)}_q$,
where $\{\varepsilon_j: 0 \leq j \leq m\}$ is the natural basis for $\mathbb
C^{m+1}$. Each $e_j^n \varepsilon_{jq}$ is a function on $\mathbb D$ taking values  in 
$\mathbb C^d$; its part in $\mathbb C^{d_j}$ is $\varepsilon_j \otimes \varepsilon^{(j)}_q$, and 
its part in every other $\mathbb C^{d_k}$ ($k\not = j)$ is $0$.  Defining
\begin{equation} \label{onbHetaY}
\mathbf e_{jq}^n := \Gamma^{(\eta, Y)}\big ( e_j^n \varepsilon_{jq}\big ), 
\end{equation}
we have an orthonormal basis for $\mathcal H^{(\eta,Y)}$.

We identify the ``$K$ -types'' in $\mathcal H^{(\eta, Y)}$,  that is, the subspaces on which the representation 
$U$ restricted to $\tilde{\mathbb K}$ acts by scalars.  For $k_\theta\in \tilde{\mathbb K}$ given by 
$k_\theta(z) = e^{i\theta} z$, we have
$D_{k_\theta}^{(\eta+j)} e^n_j = e^{-i\theta(\eta+j+n)}e^n_j$ on $\mathbb A^{(\eta+j)}$.   
By the intertwining property of $\Gamma^{(\eta,Y)}$, the basis elements of $\mathcal H^{(\eta,Y)}$ 
then satisfy $U_{k_\theta}\mathbf e_{jq}^n = 
e^{-i\theta(\eta+j+n)}\mathbf e_{jq}^n$.  It follows that the subspace
$$\mathcal H^{(\eta,Y)}(n):=\{f \in \mathcal H^{(\eta,Y)}: U_{k_\theta}f =
e^{-i\theta(\eta+n)}f\}$$
is spanned by the basis elements
$\{\mathbf e_{jq}^{n-j}: \,1 \leq q \leq d_j,\:0 \leq j \leq \min(m,n) \}$
and $\mathcal H^{(\eta,Y)}$ equals the direct sum $\oplus_{n \geq 0}\mathcal H^{(\eta,Y)}(n)$.

Clearly, the operator $M$ maps each $\mathcal H^{(\eta,Y)}(n)$ to $\mathcal H^{(\eta,Y)}(n+1)$. 
We write $M(n)$ for the matrix of the restriction of $M$ to 
$\mathcal H^{(\eta,Y)}(n)$, that is,  
\begin{equation} \label{Mblock}
M \mathbf e_{jq}^{n-j} =  \sum_{\ell, r}M(n)_{(\ell r)(j q)} \mathbf
e^{n+1-\ell}_{\ell r}.
\end{equation}
We write $e^{n-j}_{(jq),(st)}(z)$ for the $(s,t)$ - component ($0\leq s \leq \min(m,n)$, $1 \leq t \leq d_s$) 
of the function $\mathbf e_{jq}^{n-j}$ taking values in $\mathbb C^d$.  This can be regarded as a matrix of 
monomials in $z$.   The coefficients of these monomials form a numerical matrix which we denote by $E(n)$.

Applying the operator $M$, which is multiplication by $z$, to the monomials does not change their coefficients.  
Therefore, equation \eqref{Mblock} amounts to the matrix equality
\begin{equation}\label{Eblock}
E(n) = E(n+1) M(n).
\end{equation}
We use \eqref{onbHetaY} to compute $E(n)$ explicitly.  The part in $\mathbb C^{d_j}$ of the vector valued function 
$e^{n-j}_j\varepsilon_{jq}$ is  $e^{n-j}_j\varepsilon_q^{(j)}$, its part in $\mathbb C^{d_k}$ with $k\not = j$ is $0$.  
So \eqref{Gammadefn} gives, for the part of $\mathbf e^{n-j}_{jq}$ ($0\leq j \leq m$) in $\mathbb C^{d_\ell}$,
\begin{equation} \label{OrthCoeff}
\big (\mathbf e^{n-j}_{jq}(z)\big )_{\ell } =
\begin{cases}
c(\eta, \ell, j, n) z^{n-\ell} \big (Y_\ell \cdots Y_{j+1}\big ) \varepsilon_q^{(j)} z^{n-\ell}
&\mbox{if~} \ell \geq j\\
0 & \mbox{if~}\ell < j,
\end{cases}
\end{equation}
where the constant $c(\eta, \ell, j, n)$ is the coefficient of $z^{n-\ell}$ in
$$
\frac{1}{(\ell-j)!}\frac{1}{(2\eta +2j)_{\ell-j}}  \big (\frac{d}{dz} \big )^{\ell-j}e^{n-j}_j(z)
= \frac{1}{(\ell-j)!}\frac{1}{(2\eta +2j)_{\ell-j}} \sqrt{\frac{(2\eta+2j)_n}{n!}}
\big (\frac{d}{dz} \big )^{\ell-j} z^{n-j}.
$$

We can regard $E(n)$ as a block matrix with blocks $E(n)_{j \ell}$ of size $d_j \times d_\ell$.  
The $(q,r)$ entry of $E(n)_{j \ell}$ being $E(n)_{(jq)(\ell r)}$ defined above.  Then equation \eqref{OrthCoeff} says that
$$
E(n)_{j \ell } =
\begin{cases}
c(\eta, \ell, j, n)  Y_\ell \cdots Y_{j+1}
&\mbox{if~} \ell \geq j\\
0 & \mbox{if~}\ell < j.
\end{cases}
$$
Now, we consider the behavior of $c(\eta, \ell, j, n)$ for large $n$.
First, since
\begin{eqnarray*}
\sqrt{\frac{(2\eta +2j)_{n-j}}{(n-j)!}}\,\big (\frac{d}{dz} \big )^{\ell-j} z^{n-j}
&=&\frac{\sqrt{(n-j)!(2\eta+2j)_{n-j}}}{(n-\ell)!}\,z^{n-\ell},
\end{eqnarray*}
it follows that 
$$
c(\eta,\ell,j,n) = \frac{1}{(2\eta+2j)_{\ell-j}(\ell-j)!} \frac{\sqrt{\Gamma(n-j+1)\Gamma(2\eta+j+n)}}
{\sqrt{\Gamma(2\eta+2j)}\Gamma(n-\ell+1)}.
$$
From Stirling's formula,  we obtain 
\begin{eqnarray*}
c(\eta,\ell,j,n) &\sim&\frac{1}{\sqrt{\Gamma(2\eta+2j)(2\eta+2j)_{\ell-j}(\ell-j)!}}
\frac{\sqrt{(e^{-n}n^{n-j+ \frac{1}{2}})(e^{-n}n^{n+2
\eta+j-\frac{1}{2}})}}{e^{-n}n^{n-\ell+\frac{1}{2}}} \\
&\sim&\frac{\sqrt{\Gamma(2 \eta+2j)}}{\Gamma((\ell-j+1))\Gamma(2\eta+2j+\ell)}
n^{\eta-\frac{1}{2}+\ell}.
\end{eqnarray*}
If we define the block matrix $E$ by 
$$
E_{\ell j}=\begin{cases}
\frac{\sqrt{\Gamma(2 \eta+2j)}}{\Gamma((\ell-j+1))\Gamma(2\eta+2j+\ell)}
 Y_\ell\cdots Y_{j+1} & \mbox{if~}\ell\geq j\\
0 &\mbox{if~} \ell < j
\end{cases}
$$
and the diagonal block matrix $D(n)$ by $D(n)_{\ell \ell}= n^\ell I_{d_\ell}$ then we can write our result as
$$
E(n) \sim n^{\eta-\frac{1}{2}}D(n)E
$$
From \eqref{Eblock}, for large $n$, it follows that 
\begin{eqnarray}\label{HSPerturb}
M(n) &=& E(n+1)^{-1}E(n) \nonumber\\ &\sim& \big ( \frac{n}{n+1}\big
)^{\eta-\frac{1}{2}} E^{-1} D(n+1)^{-1}D(n)E \nonumber\\&\sim& I + O(1/n).
\end{eqnarray}
Therefore, the operator $M$ which is a ``weighted block shift'' is the direct sum of an ordinary 
(unweighted) block shift and a Hilbert - Schmidt operator. Hence  $M$ is bounded and standard results 
from Fredholm theory ensure that the adjoint operator $M^*$ is in the Cowen-Douglas class ${\rm B}_d(\mathbb D)$.
\end{proof}
The similarity classes of the homogeneous Cowen-Douglas operators are easily described in terms of the 
parameter $\eta$ and the multiplicities $d_0, \ldots ,d_m$. For a somewhat smaller class of operators, 
the similarity classes were described in \cite{SSR1}.  
\begin{Theorem} \label{Halmos}
The multiplication operator $M$ on $\mathcal H^{(\eta, Y)}_I$  and on $\mathcal H^{(\eta^\prime, Y^\prime)}_I$
are similar if and only if the blocks in $Y$ and $Y^\prime$ are of the same size and $\eta=\eta^\prime$.
\end{Theorem}
\begin{proof}
To prove the forward direction, first we show that $\Gamma^{(\eta, Y)}$ maps $\mathbf A^{(\eta)}$ onto
itself, that is, $\mathbf A^{(\eta)} = \mathcal H^{(\eta, Y)}_I$ as linear spaces. The derivative
$\frac{d}{dz}:\mathbb A^{(\alpha)}\to \mathbb A^{(\alpha+1)}$
defines a surjective bounded linear operator for any $\alpha >0$. For any $f \in \mathbf A^{(\eta)}$,
$$
(\Gamma^{(\eta, Y)} f)_\ell = \sum_{j=0}^\ell (\Gamma^{(\eta, Y)} f_j)_\ell
$$   
and \eqref{Gammadefn} shows that each term of the sum is in $d_\ell \mathbb A^{(\eta+\ell)}$.
On the other hand, given $g=(g_1,\ldots , g_m) \in \mathbf A^{(\eta)}$, we find $f \in \mathbf A^{(\eta)}$ 
satisfying $\Gamma^{(\eta, Y)} f =  g$.  The functions $f_0, \ldots ,f_d$ are determined recursively. 
Suppose, we have already determined $f_j$, $j < \ell$.  Then from the definition of the map  $\Gamma^{(\eta, Y)}$, 
we see that taking 
$$
f_\ell = g_\ell - \sum_{j=0}^{\ell-1} (\Gamma^{(\eta, Y)} f_j)_\ell
$$
we have the required $f$.  Clearly, $M:\mathbf A^{(\eta)} \to \mathbf A^{(\eta)}$ is  similar to $M:\mathbf
\mathcal H_I^{(\eta, Y)} \to \mathcal H_I^{(\eta, Y)}$ via the map $f \mapsto f$, which is bounded and
invertible by
the Closed graph theorem.

For the proof in the other direction, 
let $K(n)\subseteq \mathbf A^{(\eta)}=\oplus_{j=0}^m d_j\mathbb A^{(\eta+j)}$ be the linear span of the vectors 
$\{e^n_{jq}: 0\leq j \leq m,\,\, 1 \leq q \leq d_j\}$.  
The multiplication operator $M$ on $\mathbf A^{(\eta)}$ maps $K(n)$ into $K(n+1)$. If $M_n$ is the matrix
representing $M_{|\,K(n)}:K(n) \to K(n+1)$ then $M$ is a block shift with blocks $\{M_n:\, n\geq 0\}$, which are
diagonal matrices of size $d\times d$.
Let $M^\prime$ be the multiplication operator on $\mathbf A^{(\eta^\prime)} = \oplus_{j=0}^{m^\prime}
d^\prime_j \mathbb A^{(\eta^\prime + j)}$ with a similar block decomposition. Assume without loss of generality
that  $\eta^\prime > \eta$.
Suppose $L:\mathbf A^{(\eta)}\to \mathbf A^{(\eta^\prime)}$ is a bounded and invertible linear map 
consisting of $d\times d$ blocks with $LM = M^\prime L$. Then $d=d_0+\cdots + d_m = {\rm codim} ({\rm ran}\, M
)= {\rm codim} ({\rm ran}\, M^\prime ) = d_0^\prime + \cdots + d_{m^\prime}^\prime$.

It then follows that $L_{0k}=0$ for all $k \geq 1$ and consequently $L_{00}$ 
is non-singular.  We also have $L_{n\,n} M_{n-1} = M^\prime_{n-1}L_{n-1\, n-1}$ from which it follows that
$$
L_{n\,n} = M^\prime_{n-1}\cdots M^\prime_0 L_{00} M_0^{-1} \cdots M_{n-1}^{-1} = F^\prime_nL_{00}F_n^{-1},
$$
where $F_n = M_0 \cdots M_{n-1}$ and $F^\prime =  M^\prime_{n-1}\cdots M^\prime_0$ are diagonal matrices. 
The diagonal elements are 
$$F(n)_{k k} = \sqrt{\frac{(2\eta+2j(k))_n}{n!}}\,\, \Big (\mbox{\rm respectively,\,} 
F^\prime(n)_{\ell\,\ell} = \sqrt{\frac{(2\eta^\prime+2j^\prime(\ell))_n}{n!}} \Big),$$
where $j(k)=j$ if $d_0+\cdots +d_{j-1} < k \leq d_0 + \cdots +d_j$.  By Stirling's formula, we have
$$
\big (L_{nn}\big )_{\ell k}=  \big (F^\prime_n\big )_{\ell\ell}\big (L_{00}\big )_{\ell k}\big (F_n^{-1}\big )_{kk}\sim
n^{\eta^\prime-\eta + j^\prime(\ell) - j(k)}\big (L_{00}\big )_{\ell k}.
$$

Since $L_{00}$ is nonsingular, for any $k$ with $j(k)=0$, there is an $\ell$ such that $\big (L_{00}\big )_{\ell
k} \not= 0$. 
Now, unless $\eta = \eta^\prime$, we have $\big (L_{nn} \big)_{\ell k} \to \infty$ contradicting the boundeness of $L$.  
Therefore, we have $\eta = \eta^\prime$ and $\big (L_{nn} \big)_{\ell k} \sim n^{j^\prime(\ell) - j(k)}
\big (L_{00}\big )_{\ell k}$. Take all those $k$ for which $j(k)=0$.  For each of these, we can find a
different $\ell_k$  
such that $\big (L_{00}\big )_{\ell_k k} \not = 0$.  (The columns of the nonsingular matrix $L_{00}$ with these
indices 
are linearly independent and therefore cannot have only zeros in more than $d-k$ slots.) Again, unless $j^\prime(\ell_k)=0$, 
we have $\big (L{nn}\big )_{\ell_k k} \to \infty$. This shows that $d_0^\prime \geq d_0$.  Similarly, 
$d_j^\prime \geq d_{j}$, $1\leq j\leq m$.  From the equality $\sum_{j=0}^{m\prime} d_j^\prime =
\sum_{j^\prime=0}^{m^\prime} d_{j}$,  it follows that  $m^\prime = m$ and $d^\prime_j = d_j$ for $j=1, \ldots
,m$.
\end{proof}
The following Corollary, the proof of which is evident, implies that polynomially bounded homogeneous operators
in the Cowen-Douglas class are similar to contractions.
\begin{Corollary}
A homogeneous operator in the Cowen-Dougls class is either similar to a contraction or it is not power bounded. 
\end{Corollary}

\section{Examples}
In this last Section, we discuss how some formerly known examples fit into the present framework.
\subsection{\sf The case of $d_0=d_1=\cdots = d_m=1$} This case was already studied in \cite{KM}. 
Here each $Y_j$ is a number, non-zero in the irreducible case. 
The unitaries implementing the equivalence are diagonal, and clearly the conjugacy class 
$[Y]$ under these has exactly one representative with $y_j > 0$, $1 \leq j \leq m$.  
The positive $m+1$ - tuples satisfying the condition given in Theorem \ref{Pcond} give a 
parametrization of homogeneous Cowen-Douglas operators.  For each one, 
$K(0,0)=I$ and $J_g^{(\eta, Y)}$ is given by the formula \eqref{multcovering}.

Another good parametrization is possible with the aid of Theorem \ref{equiv}.  All possible 
$Y$-s are now conjugate under diagonal unitaries $A$, so we may fix an arbitrary 
$Y^{(0)}$ (for example, $y_j =1$ for all $j$, or, as in \cite{KM}, $y_j=j$ for all $j$).  
Take any positive diagonal matrix $N$ with diagonal elements $1=\mu_0, \mu_1,\ldots, \mu_m$.  
By Proposition \ref{prop2.2}, $Y^{(0)}, N$ and $Y^{(0)}, N^\prime$ 
give isomorphic vector bundles if and only if $A$ is diagonal and hence $N=N^\prime$.  
It follows that the positive numbers $\eta, \mu_1,\ldots, \mu_m$
give a parametrization of the homogeneous operators in the Cowen-Douglas class 
$B_{m+1}(\mathbb D)$.  Here $J_g^{(\eta, Y^{(0)})}$ depends only on $\eta$ and
$K_N^{(\eta, Y^{(0)})}(0,0)$ is given by the formula in Theorem \ref{equiv}.  
This is the parametrization used in \cite{KM}.

In the case $m=1$, for any  $d_0$ and $d_1$,  
the class $[Y]$ always contains a member for which $Y$ is diagonal.  So, the corresponding bundle 
is reducible unless $d_0=d_1=1$.  When $m=2$, it is easy to see that $d_0=2$ or $d_2=2$ gives only 
reducible bundles.  So, the first non-trivial case occurs (apart from the case $d_0=d_1=d_2=1$, 
which has been dealt with previously) when $d_0=d_2=1,\,d_1=2$.

\subsection{\sf The case of $(d_0,d_1,d_2) = (1,2,1)$} For this case, again there are two natural parametrizations.
Conjugating $Y$ with a block-diagonal unitary having blocks $u_0,U_1,u_2$ changes 
$Y_1,Y_2$ into $U_1Y_1u_0^{-1}$, $u_2Y_2U^{-1}_1$. Now, $U_1$ can be chosen so that 
$Y_1=\Big ( \begin{smallmatrix} a \\ \\0 \end{smallmatrix} \Big )$. Then $u_0,u_2$ and a 
scalar factor in front of $U_1$ can be found with $a \geq 0$ and 
$Y=\big ( \begin{matrix} b & c \end{matrix} \big )$ with $b, c \geq 0$.  
We   have irreducibility if and only if $a,b,c \not = 0$ and no two such triples give equivalent $Y$-s.  
So, we have a parametrization of the irreducible $E^{(\eta, Y)}$ by four arbitrary non-zero parameters.  
There is a reproducing kernel (and hence an operator in $B_4(\mathbb D)$) if and only if the the right hand side of  
the equation \eqref{Betafn} is positive; in terms of the parameters, this is
\begin{eqnarray*}
a^2 &<& 2 \eta\\
b^2 &<& \frac{2\eta+2}{1-\frac{a^2}{2(2\eta+1)}}\\
c^2 &<& 2\eta+2
\end{eqnarray*}
The positive quadruple $(\eta,a,b,c)$ subject to this condition parametrizes the homogeneous operators in $B_4(\mathbb D)$.  
In each case, $K(0,0)=I$ and $J_g$ can be expressed in terms of the parameters using \eqref{multcovering}.

The other parametrization of the $(d_0,d_1,d_2) = (1,2,1)$ case is found using Theorem \ref{equiv}.  
Simple arguments show that $Y$ can always be conjugated by a block diagonal $A$ so that  
$Y_1 = \Big ( \begin{smallmatrix} 1 \\ \\0 \end{smallmatrix} \Big )$ and $Y_2 = \big ( \begin{matrix} 1 & 0 \end{matrix} \big )$ 
or $\big ( \begin{matrix} 0 & 1 \end{matrix} \big )$.  When $Y_2 = \big ( \begin{matrix} 0 & 1 \end{matrix} \big )$, 
the bundle will be reducible for any choice of Hermitian structure. So, we can fix $Y^{(0)}$ with 
$Y_1 = \Big ( \begin{smallmatrix} 1 \\ \\0 \end{smallmatrix} \Big )$ and $Y_2 = \big ( \begin{matrix} 1 & 0 \end{matrix} \big )$. 
The block diagonal $A$ that conjugates this $Y^{(0)}$ to itself is a diagonal matrix with $(p,p,q,p)$ on the diagonal.  
If $N$ is any positive diagonal $\mbox{\rm diag}\,(n_0,N_1,n_2)$ with $n_0=1$, $N_1 = 
\Big ( \begin{matrix} \alpha & \beta\\ \bar{\beta} & \gamma \end{matrix} \Big )$ and $n_2 \geq 0$, then we can ensure 
$n_1=1=n_2$ and $\alpha,\,\beta,\,\gamma > 0$ after conjugating by an $A$.  
Thus the homogeneous bundles with reproducing kernel (and hence the homogeneous operators in $B_4(\mathbb D)$) of type $(1,2,1)$ 
are now parametrized by four positive numbers $(\eta, \alpha,\beta,\gamma)$ subject to the condition $\beta^2 < \alpha \gamma$.

By a different construction, a large subset of these examples already occurs in \cite{SSR}.

\bibliographystyle{amsplain}
\bibliography{homker}

\end{document}